\documentclass[12pt]{amsart}

\usepackage{amssymb}
\usepackage{amsmath}
\usepackage{amscd}
\usepackage{float}
\usepackage{amsfonts}
\usepackage{pb-diagram}
\usepackage{eufrak}
\usepackage[all, poly, knot]{xy}
\usepackage{graphicx, graphics}

\input xy
\xyoption{arc}
\xyoption{all}
\xyoption{knot}

\makeatletter
\long\def\M#1{\leavevmode\setbox\@tempboxa\hbox{#1}\@tempdima\fboxrule
    \advance\@tempdima \fboxsep \advance\@tempdima \dp\@tempboxa
   \hbox{\lower \@tempdima\hbox
  {\vbox{\hrule \@height \fboxrule
          \hbox{  \hskip\fboxsep
          \vbox{\vskip\fboxsep \box\@tempboxa\vskip\fboxsep}\hskip
                 \fboxsep\vrule \@width \fboxrule}%
                  }}}}
\makeatother

\let \ttorg \tt \def \tt{\ttorg \obeyspaces}

\begin{document}

\title[Tangle Insertion Invariants]{Tangle Insertion Invariants for Pseudoknots, Singular Knots, and Rigid Vertex Spatial Graphs}

\author{Allison K. Henrich}
\address{ Department of Mathematics, Seattle University,
 901 12th Ave. Seattle, WA 98122-1090, U.S.A.}
\email{henricha@seattleu.edu}
\urladdr{http://fac-staff.seattleu.edu/henricha}

\author{Louis H. Kauffman}
\address{ Department of Mathematics, Statistics and
 Computer Science, University of Illinois at Chicago,
 851 South Morgan St., Chicago IL 60607-7045, U.S.A.}
\email{kauffman@math.uic.edu}
\urladdr{http://www.math.uic.edu/$\sim$kauffman/}




\date{}
\maketitle

\thispagestyle{empty}

\begin{abstract}
{\em The notion of a pseudoknot is defined as an equivalence class of knot diagrams that may be missing some crossing information. We provide here a topological invariant schema for pseudoknots and their relatives, 4-valent rigid vertex spatial graphs and singular knots, that is obtained by replacing unknown crossings or vertices by rational tangles. }
\end{abstract}

\section{Introduction}

\begin{figure}[htbp]
\begin{center}
\includegraphics[scale=.75]{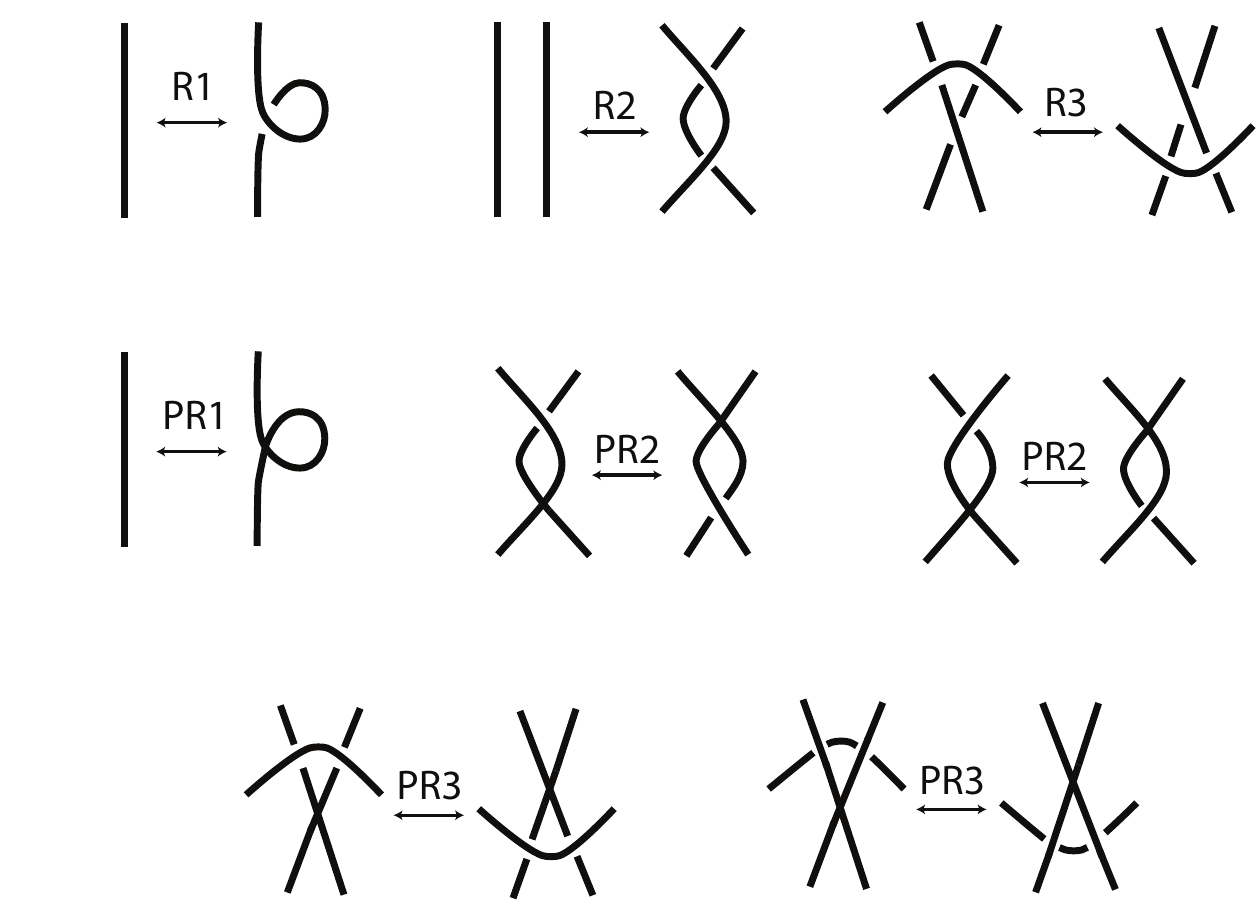}
\caption{Classical and pseudo-Reidemeister moves}
\label{Rmoves}
\end{center}
\end{figure}

Pseudoknots and pseudolinks are knots and links about which we have incomplete information. This incompleteness is expressed in diagrams by the appearance of {\it precrossings} that have no over or under designation. In~\cite{pseudo}, pseudoknots are defined as equivalence classes of knot or link diagrams (called {\em pseudodiagrams} in~\cite{hanaki}) where some crossing information may be missing. We picture the precrossings that lack definitive over/under information, as undecorated self-intersections. These unknown crossings can be involved in Reidemeister-like moves in predictable ways. (See Figure~\ref{Rmoves}.)

In terms of their behavior, pseudoknots act like 4-valent rigid vertex spatial graphs (as in~\cite{lou_graphs}) or singular knots (as in~\cite{singular}). The one important difference between pseduoknots and these other objects is entirely characterized by the PR1 move. The PR1 move allows us to eliminate a single unknown crossing at a curl. The reason this move is reasonable for pseudoknots is that, regardless of whether we replace the precrossing with a positive or negative classical crossing, it can be removed with an R1 move. On the other hand, vertices or singularities are not removable in rigid vertex spatial graphs or singular knots.

While several new pseduoknot invariants have been introduced~\cite{gauss, heather, pseudo, color}, we introduce a schema for a large collection of new invariants. These invariants of pseudoknots, defined by replacing precrossings with tangles, can be adapted to serve as invariants for 4-valent rigid vertex spatial graphs and singular knots as well.

\section{Tangle Insertion Invariants}

In~\cite{singular}, Vassiliev showed how to extend an invariant $\mathcal{I}$ of links to the class of singular links. He did so by defining the {\em derivative} of an invariant as follows: 

\smallskip
\begin{center}{\huge
$\mathcal{I}'($\includegraphics[height=.35in]{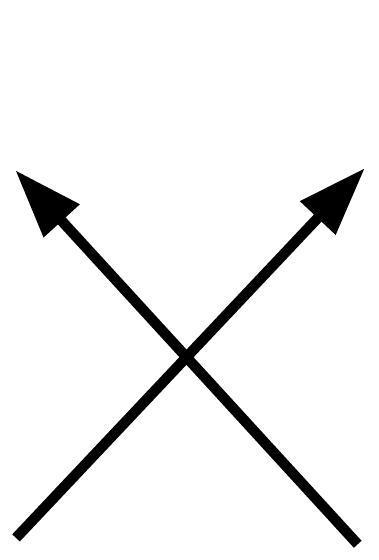}$)=\mathcal{I}($\includegraphics[height=.35in]{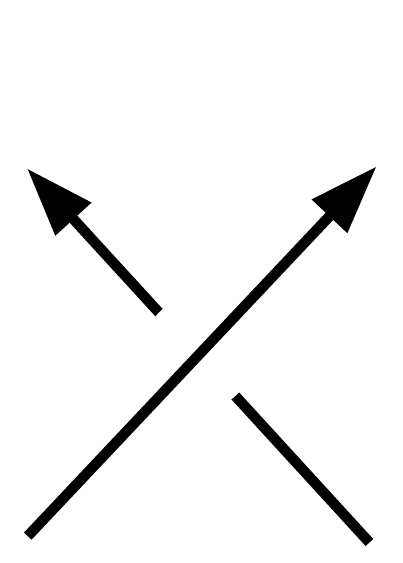}$)-\mathcal{I}($\includegraphics[height=.35in]{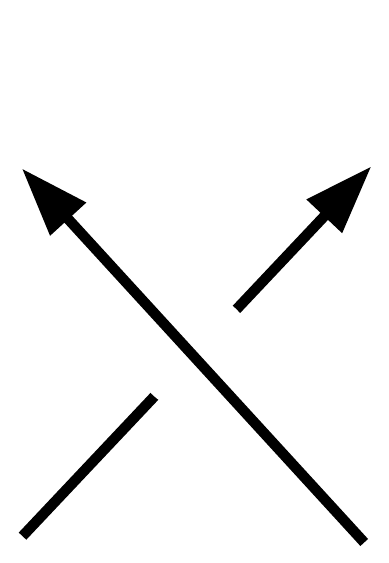}$)$}
\end{center}
\smallskip

The derivative of an invariant is an invariant of singular links with a single double-point. In general the $n$th derivative, defined recursively as follows, is an invariant of a singular link with $n$ double-points. Note that when $n=1$, $\mathcal{I}^{(n)}$ is the first derivative $\mathcal{I}'$, and $\mathcal{I}^{(n-1)}=\mathcal{I}$.

\smallskip
\begin{center}{\huge
$\mathcal{I}^{(n)}($\includegraphics[height=.35in]{Precrossing.pdf}$)=\mathcal{I}^{(n-1)}($\includegraphics[height=.35in]{Pos.pdf}$)-\mathcal{I}^{(n-1)}($\includegraphics[height=.35in]{Neg.pdf}$)$}
\end{center}
\smallskip

This idea of extending invariants by taking linear combinations of invariant values for diagrams that are related by tangle replacements can be generalized (see \cite{KaufVogel}). For instance, we could define $\overline{\mathcal{I}}$ recursively as follows.

\smallskip
\begin{center}{\huge
$\overline{\mathcal{I}}($\includegraphics[height=.35in]{Precrossing.pdf}$)=a\overline{\mathcal{I}}($\includegraphics[height=.35in]{Pos.pdf}$)+b\overline{\mathcal{I}}($\includegraphics[height=.35in]{Neg.pdf}$)+c\overline{\mathcal{I}}($\includegraphics[height=.35in]{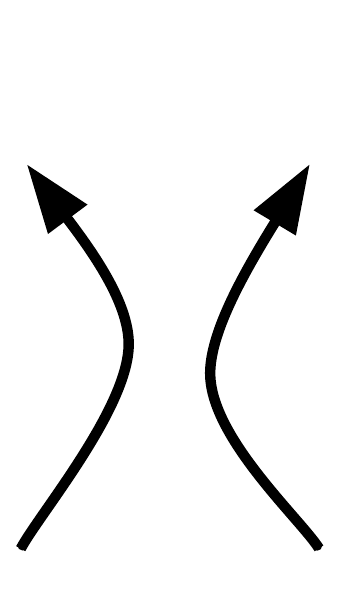}$)$}
\end{center}
\smallskip

If we are careful about which invariants $\mathcal{I}$ we use and how we choose our coefficients $a$, $b$, and $c$, this equation can be used to define an invariant of singular links, rigid vertex spatial graphs, or pseudoknots and links. Note that any choice of coefficients that yields an invariant of singular links is also an invariant of pseudoknots if it satisfies the following additional relation.

\smallskip
\begin{center}\includegraphics[height=.5in]{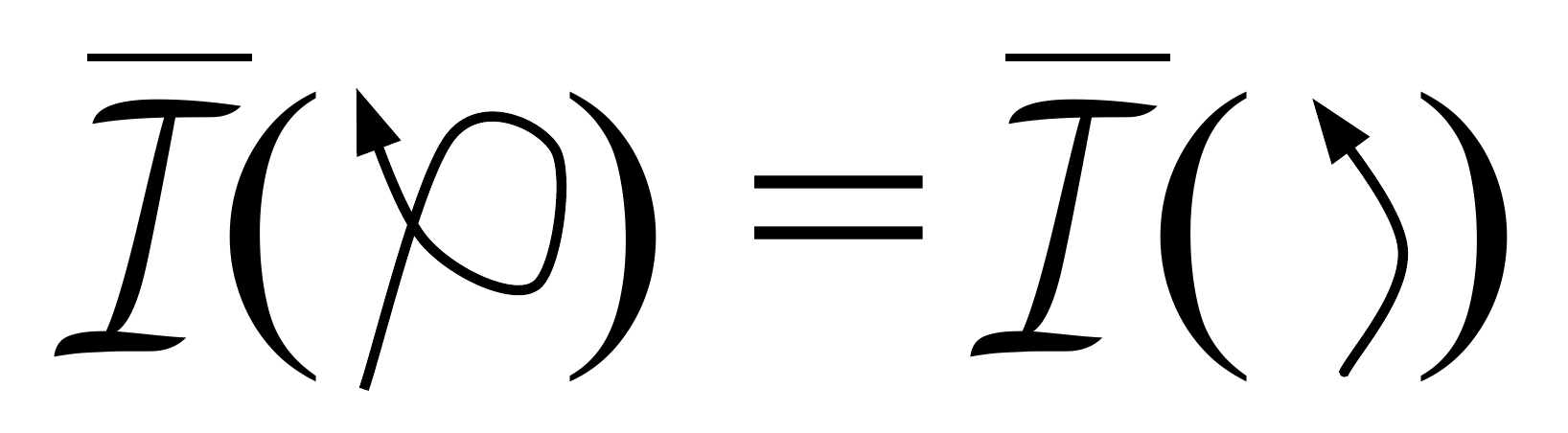}
\end{center}
\smallskip



In general, suppose $\mathcal{I}$ is a polynomial or integer, real, or complex-valued invariant of knots and links that behaves nicely under and connected sum. In particular, suppose $\mathcal{I}$ obeys the following property. 

\begin{align*}
\mathcal{I}(K\# L)&=\mathcal{I}(K)\mathcal{I}(L)
\end{align*}

\

Here, $K$ and $L$ denote arbitrary pseudoknot or link diagrams, and $\#$ indicates a connected sum. We observe that the Jones and Alexander polynomials satisfy this property, among other commonly used invariants. 

\begin{figure}[htbp]
\begin{center}
\includegraphics[scale=.5]{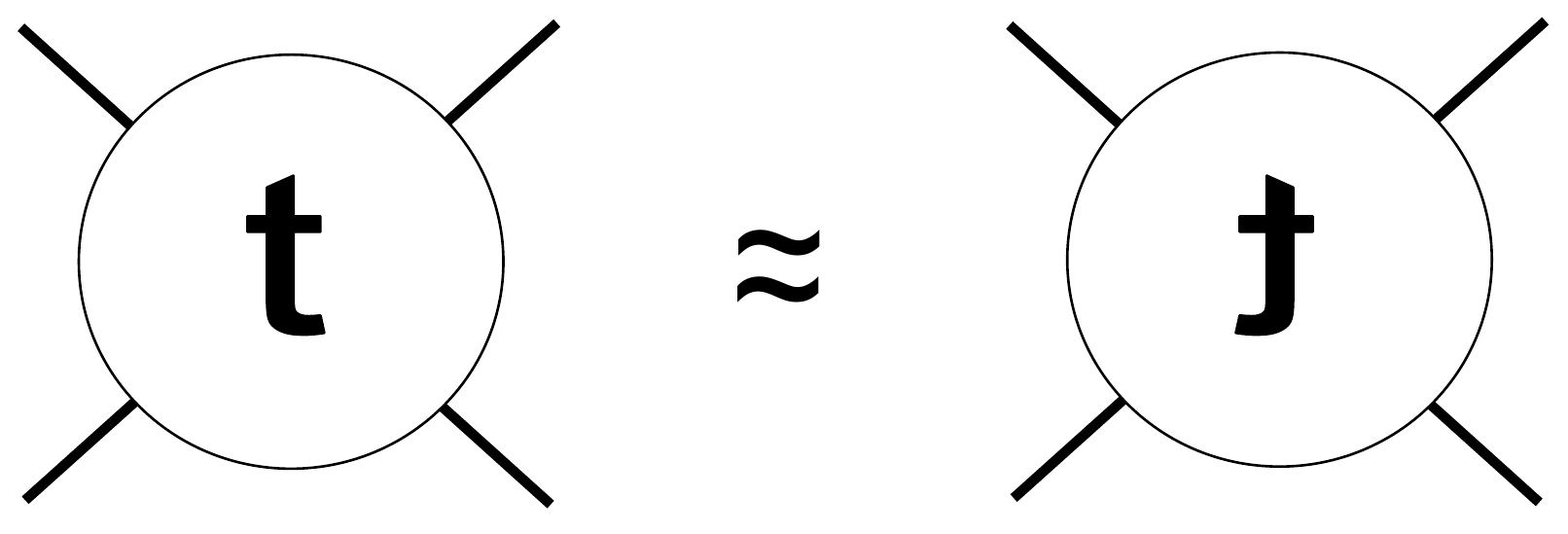}\\
{\bf Vertical Flip}\\
\bigskip

\includegraphics[scale=.5]{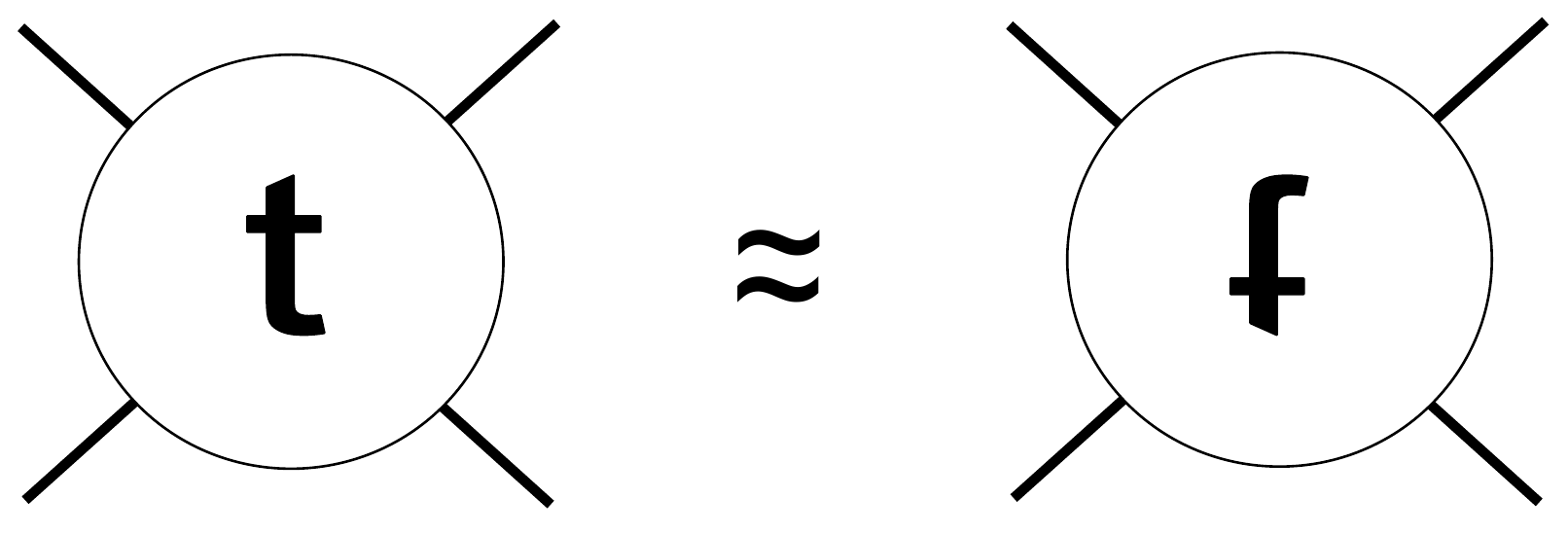}\\
{\bf Horizontal Flip}
\caption{Rational tangle symmetries}
\label{tangle_sym_free}
\end{center}
\end{figure}

\begin{figure}[htbp]
\begin{center}
\includegraphics[height=.9in]{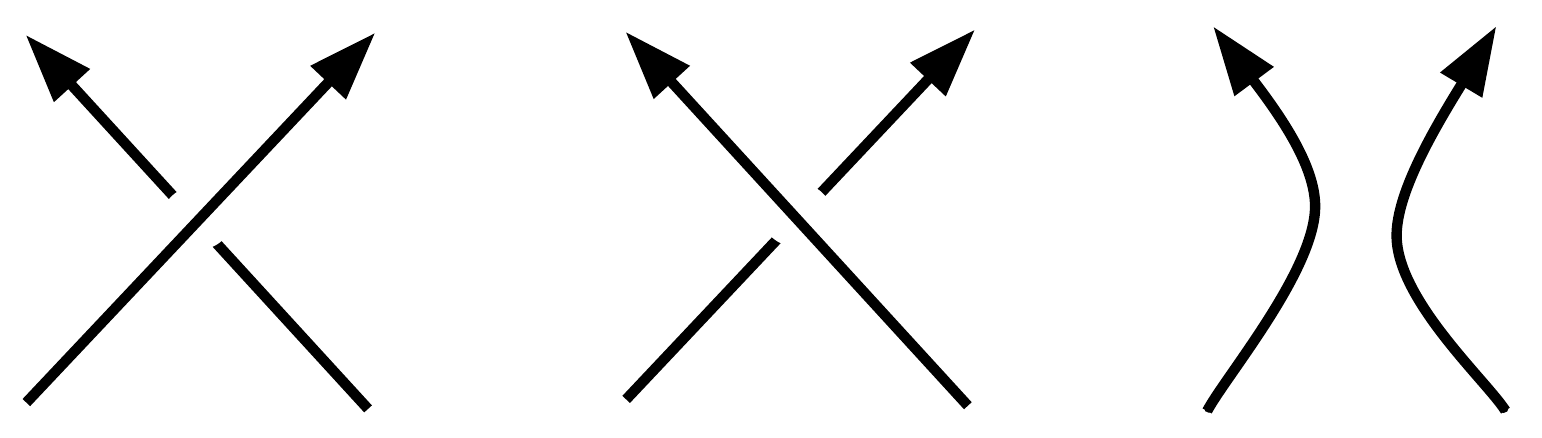}\\
(a)\hspace{1in}(b)\hspace{.9in}(c)\\
\includegraphics[height=1.5in]{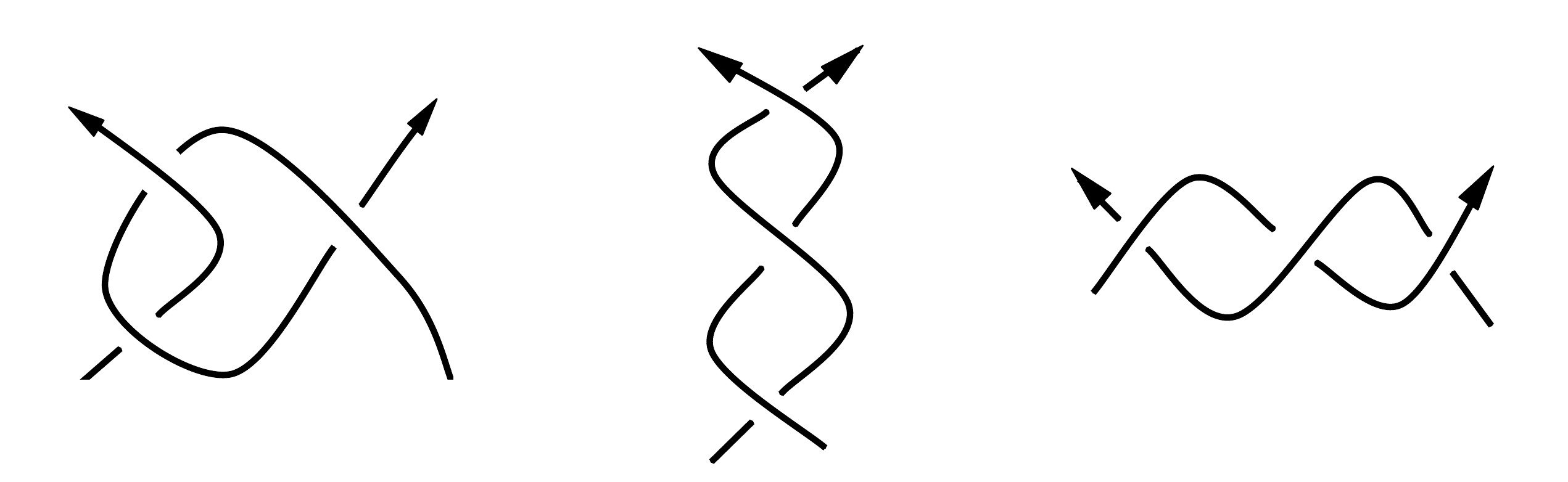}\\
(d)\hspace{1.4in}(e)\hspace{1.4in}(f)
\caption{Examples of tangles with the required symmetries}
\label{example_t}
\end{center}
\end{figure}

Once a suitable link invariant has been chosen, we choose a set of oriented, rational (2,2)-tangles, $T=\{t_1,t_2,...,t_n\}$. Recall that, since each $t\in T$ is a rational tangle, $t$ satisfies the symmetries shown in Figure~\ref{tangle_sym_free}, by the Flip Theorem for rational tangles~\cite{rational, lou_sofia1, lou_sofia2}. Some examples of oriented rational tangles are shown in Figure~\ref{example_t}.

Now that we have our desired link invariant and set of tangles, we define $\widehat{\mathcal{I}}(P)$ recursively as follows for an oriented pseudodiagram $P_c$ containing precrossing $c$.

$$\widehat{\mathcal{I}}(P_c)=\sum_{i=1}^n A_i\widehat{\mathcal{I}}(P_{c=t_i})$$

In this definition, $P_{c=t_i}$ denotes the pseudodiagram where tangle $t_i$ is inserted at precrossing $c$ (respecting the orientations of the tangles and the pseudodiagram), and the symbol $A_i$ is a variable. If a given pseudodiagram $P$ contains no precrossings, we define $\widehat{\mathcal{I}}(P)=\mathcal{I}(P)$.

The restrictions on our allowable invariants and tangles guarantee that $\widehat{\mathcal{I}}$ is an invariant of singular links and 4-valent rigid vertex spatial graphs. To guarantee invariance under {\em all} pseudo-Reidemeister moves, we need to further impose a relation among our coefficients. 

$$\overline{\mathcal{I}}(P)=\widehat{\mathcal{I}}(P)/\left<\sum_{i=1}^n A_i\mathcal{I}(D(t_i))=1\right>$$

We take $\widehat{\mathcal{I}}(P)$ modulo a linear combination of the values of the link invariant of the denominator closures of each of our tangles to ensure that $\overline{\mathcal{I}}(P)$ is invariant under PR1, and hence, is an invariant of oriented pseudoknots and links. Note that we are implicitly assuming that the value of our chosen invariant is nonzero on the denominator closure of at least one of our tangles $t_i\in T$.

\begin{figure}[htbp]
\begin{center}
\includegraphics[height=2in]{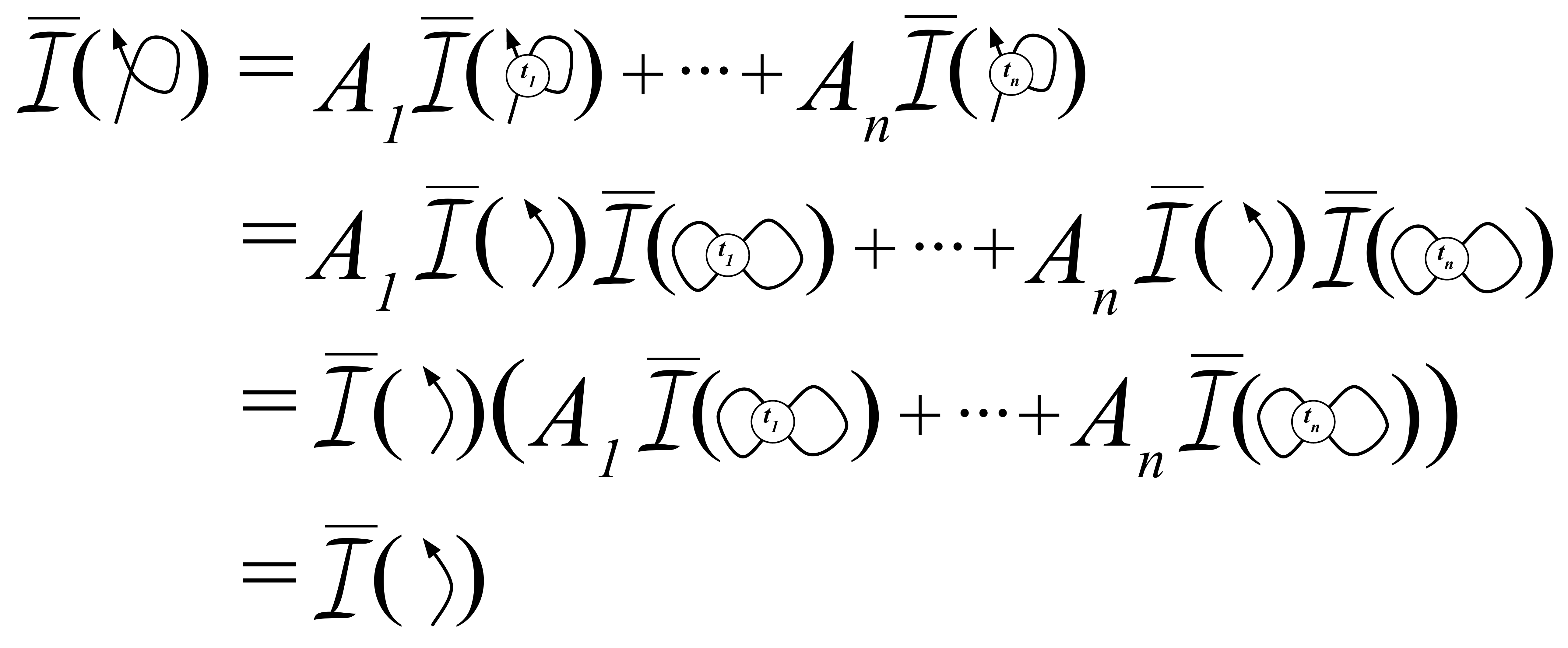}
\caption{PR1 Invariance}
\label{PR1_pf}
\end{center}
\end{figure}

To gain an understanding of why each of the restrictions imposed above is necessary, let us prove that $\overline{\mathcal{I}}(P)$ is an invariant of oriented pseudolinks. We begin by acknowledging that, since $\mathcal{I}$ is a link invariant, invariance of $\overline{\mathcal{I}}$ for classical Reidemeister moves is immediate. PR1 invariance is shown in Figure~\ref{PR1_pf}. Line 1 follows from the definition of $\overline{\mathcal{I}}$, line 2 follows from our required connected sum property of the link invariant $\mathcal{I}$ (which ensures that $\overline{\mathcal{I}}$ behaves similarly under connected sum), line 3 is an algebraic distribution, and line 4 is a consequence of the relation $\sum_{i=1}^n A_i\mathcal{I}(D(t_i))=1$. 

\begin{figure}[htbp]
\begin{center}
\includegraphics[height=3in]{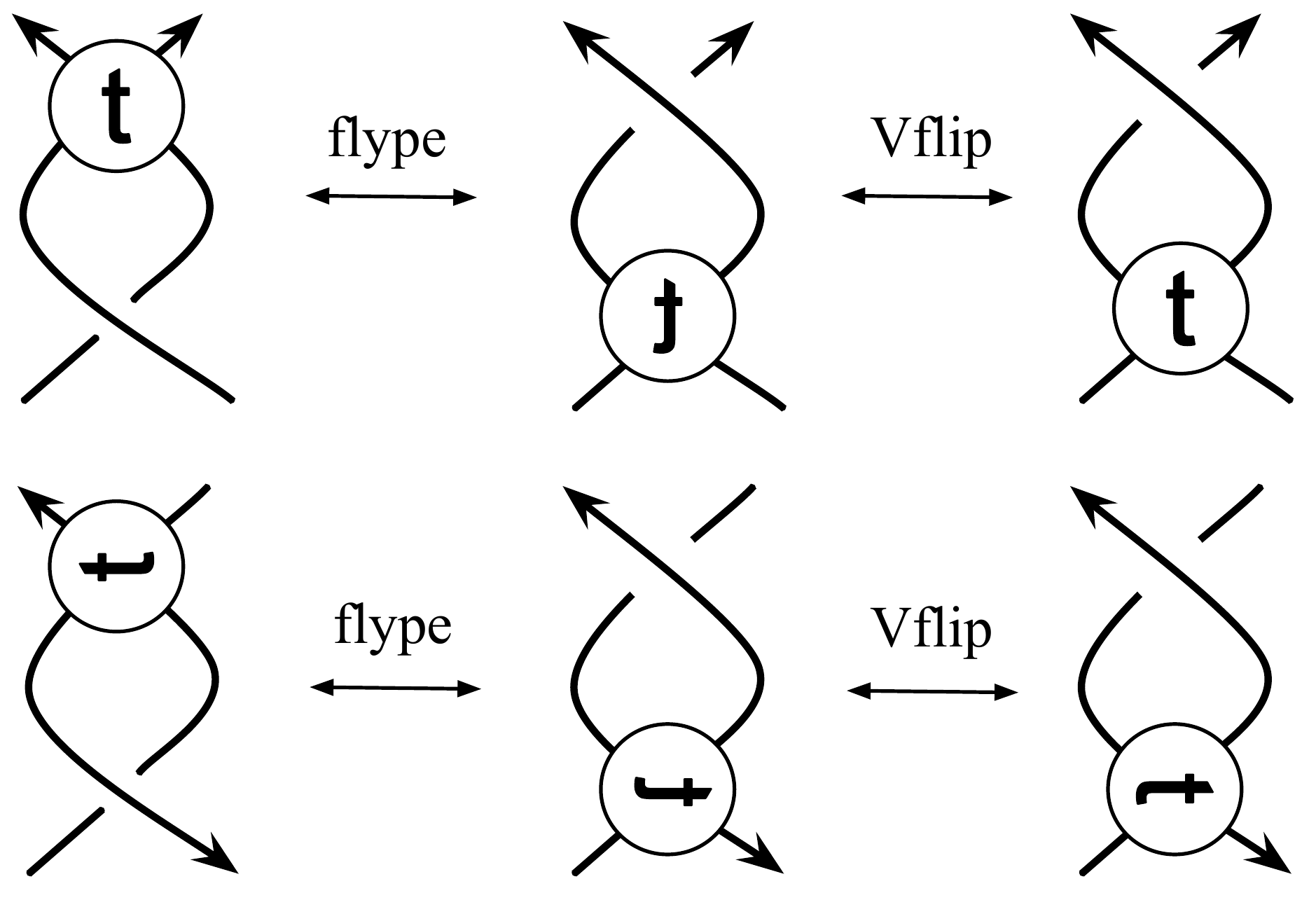}
\caption{PR2 Invariance}
\label{PR2_pf}
\end{center}
\end{figure}

The proof of invariance under the PR2 move illustrates why our tangles $t\in T$ are chosen to have the symmetries of rational tangles. A flype together with a rational tangle symmetry is required to show invariance of $\overline{\mathcal{I}}$ under each oriented PR2 move, pictured in Figure~\ref{PR2_pf}.

\begin{figure}[htbp]
\begin{center}
\includegraphics[height=1.5in]{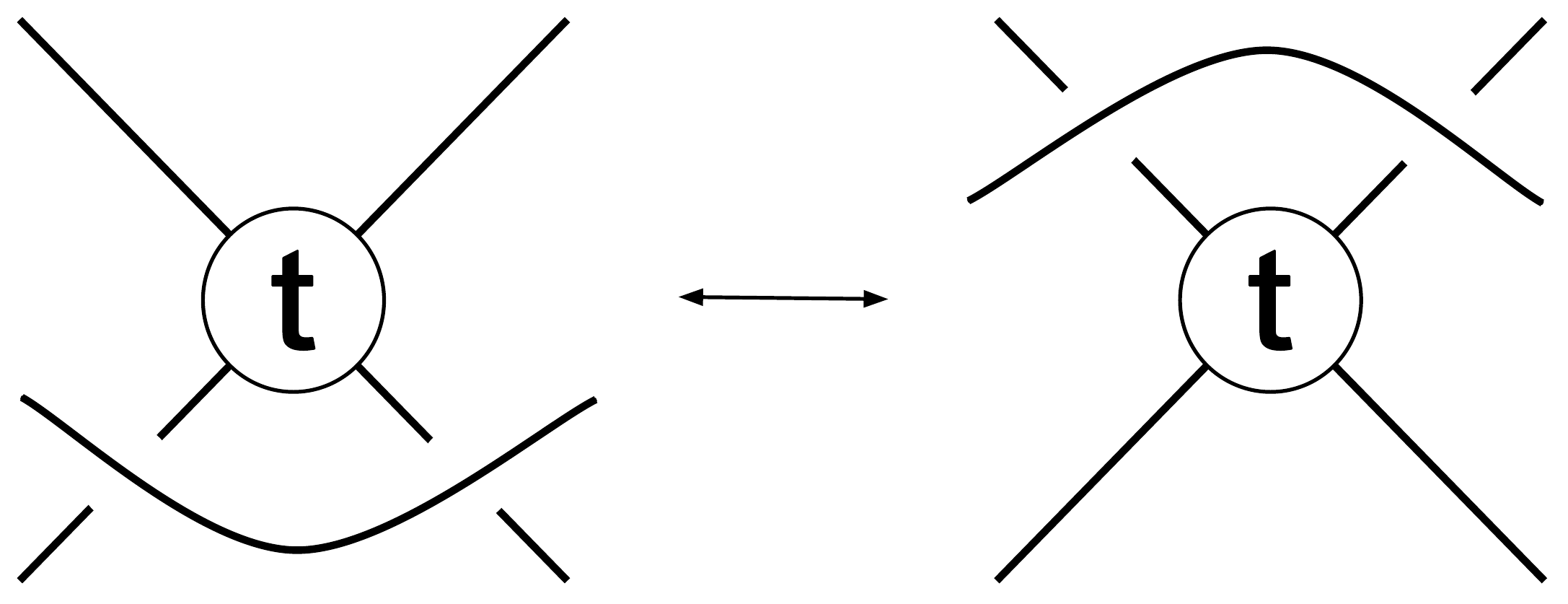}
\caption{PR3 Invariance}
\label{PR3_pf}
\end{center}
\end{figure}

Finally, PR3 invariance is straightforward. Figure~\ref{PR3_pf} illustrates how $\overline{\mathcal{I}}$ is invariant under PR3 since a strand that lies entirely above or entirely below a tangle can be moved freely past the tangle.

This completes our proof that $\overline{\mathcal{I}}$ is an invariant of pseudoknots. Note that we also proved that $\widehat{\mathcal{I}}$ is a singular link and 4-valent rigid vertex spatial graph invariant.\\

\noindent {\bf Remark.} This tangle insertion method can be generalized beyond rational tangle insertion by asking that the tangles satisfy the symmetries
shown in Figure~\ref{tangle_sym_free}.
We shall take up this aspect of the construction in another paper.\\

One obvious question to ask now that we have created a schema for creating pseudoknot invariants is: how does this schema relate to known pseudoknot invariants? One of the simplest, yet most powerful invariants of pseudoknots is the {\em weighted resolution set}, or {\em were-set}, introduced in~\cite{pseudo}. The were-set of a pseudodiagram is the set $\mathcal{S}$ of pairs $$\mathcal{S}=\{(K_1,p_1), (K_2, p_2), ..., (K_n, p_n)\}$$ of knot types $K_i$ that can be realized by some choice of crossing information for the diagram's precrossings. The number $p_i$ is the probability that knot $K_i$ will be produced if crossing information is randomly chosen, where positive and negative crossings are both equally likely. It was proven that the were-set is indeed a pseudoknot invariant. How does this invariant relate to our schema?

First, we notice that choosing crossing information for a precrossing is equivalent to inserting a basic +1 or -1 rational tangle (i.e. tangles (a) and (b) in Figure~\ref{example_t}), so let our tangle set $T$ consist of these two tangles. Next, let $\mathcal{I}$ be a knot invariant such that $\mathcal{I}(U)=1$ if $U$ is the unknot (such as the Jones polynomial or the Alexander polynomial). Then choosing the coefficients $A_1=A_2=\frac{1}{2}$ will satisfy the relation $\sum_{i=1}^n A_i\mathcal{I}(D(t_i))=1$, since this amounts to 

{\large $$
\frac{1}{2} \mathcal{I}(\xygraph{
!{0;/r1.0pc/:}
[u(.5)]!{\hcross}
!{\hcap[1]}
[l]!{\hcap[-1]}
}) 
+\frac{1}{2}\mathcal{I}(\xygraph{
!{0;/r1.0pc/:}
[u(.5)]!{\htwist}
!{\hcap[1]}
[l]!{\hcap[-1]}
})
=\frac{1}{2}+\frac{1}{2}=1,$$}

\noindent so an invariant $\overline{\mathcal{I}}$ is determined by our ingredients, following the recipe above. The invariant we have just created with our schema is equivalent to the following composition of the were-set with $\mathcal{I}$.

$$p_1\mathcal{I}(K_1)+p_2\mathcal{I}(K_2)+\cdots+p_n\mathcal{I}(K_n)$$
\smallskip

Note that, if we make a different coefficient choice, we no longer recover the were-set. For instance, making the choice $A_1=\frac{3}{4}$ and $A_2=\frac{1}{4}$ would correspond to a distinct variant of the were-set where positive crossings are chosen with probability $\frac{3}{4}$ and negative crossings are chosen with probability $\frac{1}{4}$.

\section{Examples}

Let us return to our original motivation for tangle insertion invariants: the derivative of an invariant. In our new framework, the tangles $t_1, t_2\in T$ that are used to define the derivative are tangles (a) and (b) in Figure~\ref{example_t}, respectively. Their denominator closures are both the unknot, $U$. The coefficient $A_1$ is 1 and $A_2=-1$. Notice, then, that our imposed relation $\sum_{i=1}^n A_i\mathcal{I}(D(t_i))=1$ states that $\mathcal{I}(U)-\mathcal{I}(U)=1$. But clearly $\mathcal{I}(U)-\mathcal{I}(U)=0$. So, oddly enough, our original motivating example is not an instance of our scheme after all. This is one of the reasons why the generalization is so valuable. It allows for the creation of many new invariants of pseudoknots just as our old singular knot invariants fail to be useful. 

Instead, let us consider the other example we mentioned above, namely:

\smallskip
\begin{center}{\huge
$\overline{\mathcal{I}}($\includegraphics[height=.35in]{Precrossing.pdf}$)=a\overline{\mathcal{I}}($\includegraphics[height=.35in]{Pos.pdf}$)+b\overline{\mathcal{I}}($\includegraphics[height=.35in]{Neg.pdf}$)+c\overline{\mathcal{I}}($\includegraphics[height=.35in]{Smooth.pdf}$)$}
\end{center}
\smallskip

Here, we take tangles (a), (b), and (c) from Figure~\ref{example_t} to form the tangle set $T$, and we don't yet specify values for our coefficients. 

\vspace{-.3in}

{\Huge $$
a \overline{\mathcal{I}}(\xygraph{
!{0;/r1.0pc/:}
[u(.5)]!{\hcross}
!{\hcap[1]}
[l]!{\hcap[-1]}
}) 
+b\overline{\mathcal{I}}(\xygraph{
!{0;/r1.0pc/:}
[u(.5)]!{\htwist}
!{\hcap[1]}
[l]!{\hcap[-1]}
})
+c\overline{\mathcal{I}}(\xygraph{
!{0;/r1.0pc/:}
[u(0.5)r]!{\hcap[1]}
!{\hcap[-1]}
[ll]!{\hcap[1]}
!{\hcap[-1]}
})=1$$}

\begin{figure}[htbp]
\begin{center}
\includegraphics[height=1in]{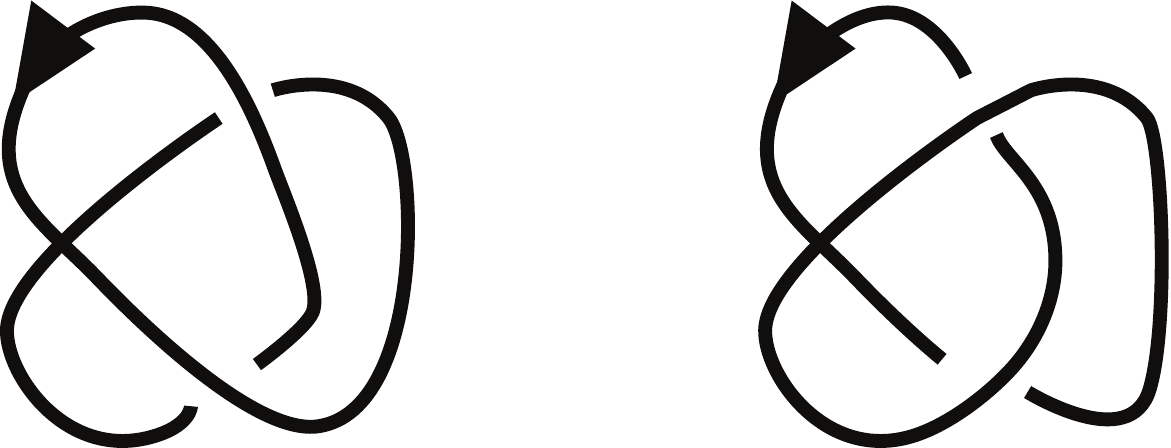}
\caption{A pseudo-trefoil and its mirror image.}
\label{tref}
\end{center}
\end{figure}

Suppose we choose our link invariant $\mathcal{I}$ to be the Jones polynomial, $J(L)$. We will compute  $J(L)$ using the bracket polynomial, $<L>$, via the relation $$J(L)=(-A^3)^{-w(L)}<L>$$ where $w(L)$ is the writhe of the link $L$ as in~\cite{K}, and $A$ is a variable. Then any choice for $a$, $b$, and $c$ satisfying $a+b+\delta c=1$ will do, where $\delta=-A^2-A^{-2}$ is the polynomial such that $\mathcal{I}(K\sqcup U)=\delta\mathcal{I}(K)$. In particular, we could choose $a=b=0$ and $c=\frac{1}{-A^2-A^{-2}}$. Note that this choice is equivalent to simply starting with the singleton tangle set consisting of the 0-tangle, (c). Let us apply this invariant to a pseudo-trefoil and its mirror image, shown in Figure~\ref{tref}.

\begin{figure}[htbp]
\begin{center}
\includegraphics[height=1in]{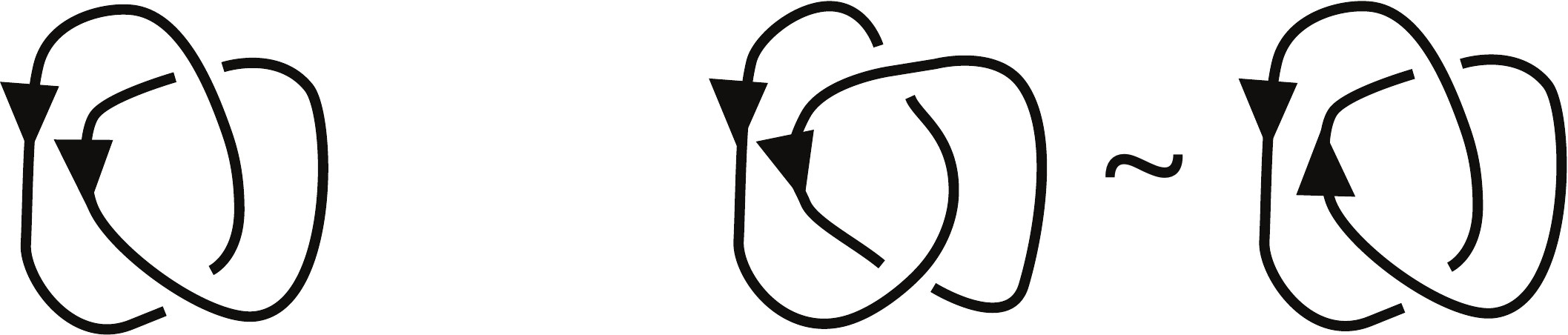}
\caption{Insertion of tangle (c) into a pseudo-trefoil and its mirror image.}
\label{tref_insert}
\end{center}
\end{figure}

When we insert tangle (c) into both pseudodiagrams, we get the links shown in Figure~\ref{tref_insert}.

The link on the left is L2a1$\{1\}$ and the link on the right is L2a1$\{0\}$. Both links have the same bracket polynomial value, $-A^4-A^{-4}$, since they are the same as unoriented links. (See~\cite{adams}, p. 151 for this computation.) But the writhe of L2a1$\{1\}$ is 2 while the writhe of L2a1$\{0\}$ is -2. Hence, the Jones polynomial of L2a1$\{1\}$ is $$(-A^4-A^{-4})(-A^3)^{-2}=-A^{-2}-A^{-10}$$ while the Jones polynomial of L2a1$\{0\}$ is $$(-A^4-A^{-4})(-A^3)^{2}=-A^{10}-A^{2}.$$ Thus, the two pseudoknots shown in Figure~\ref{tref} are distinct, since their $\overline{\mathcal{I}}$ values are $(\frac{1}{-A^2-A^{-2}})(-A^{-2}-A^{-10})$ and $(\frac{1}{-A^2-A^{-2}})(-A^{10}-A^{2})$, respectively.

\begin{figure}[htbp]
\begin{center}
\includegraphics[height=1.3in]{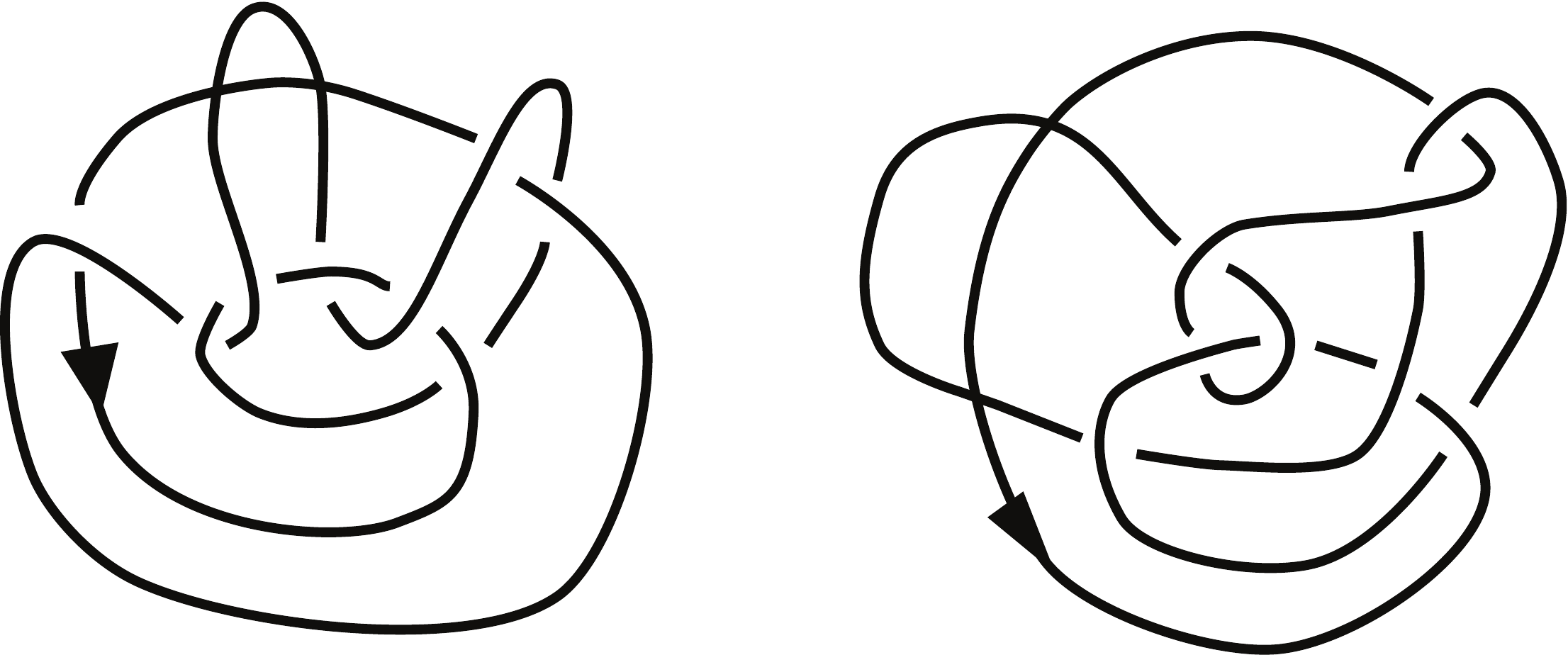}
\caption{Pseudodiagrams related to the Perko pair.}
\label{perko}
\end{center}
\end{figure}

\begin{figure}[htbp]
\begin{center}
\includegraphics[height=1.3in]{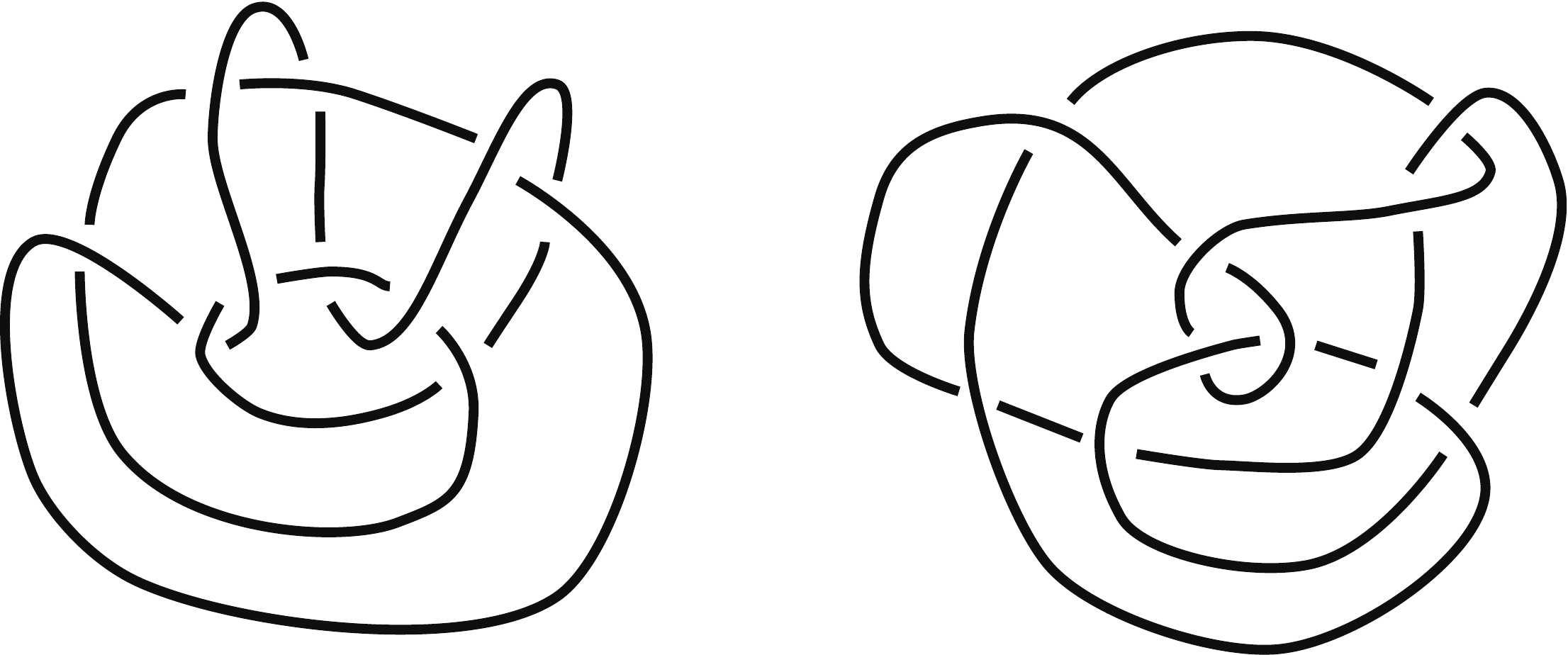}
\caption{The Perko pair.}
\label{perko_orig}
\end{center}
\end{figure}

Another interesting pair of examples, pictured in Figure~\ref{perko}, is derived from the famous Perko pair, pictured in Figure~\ref{perko_orig}. 

\begin{figure}[htbp]
\begin{center}
\includegraphics[height=1.3in]{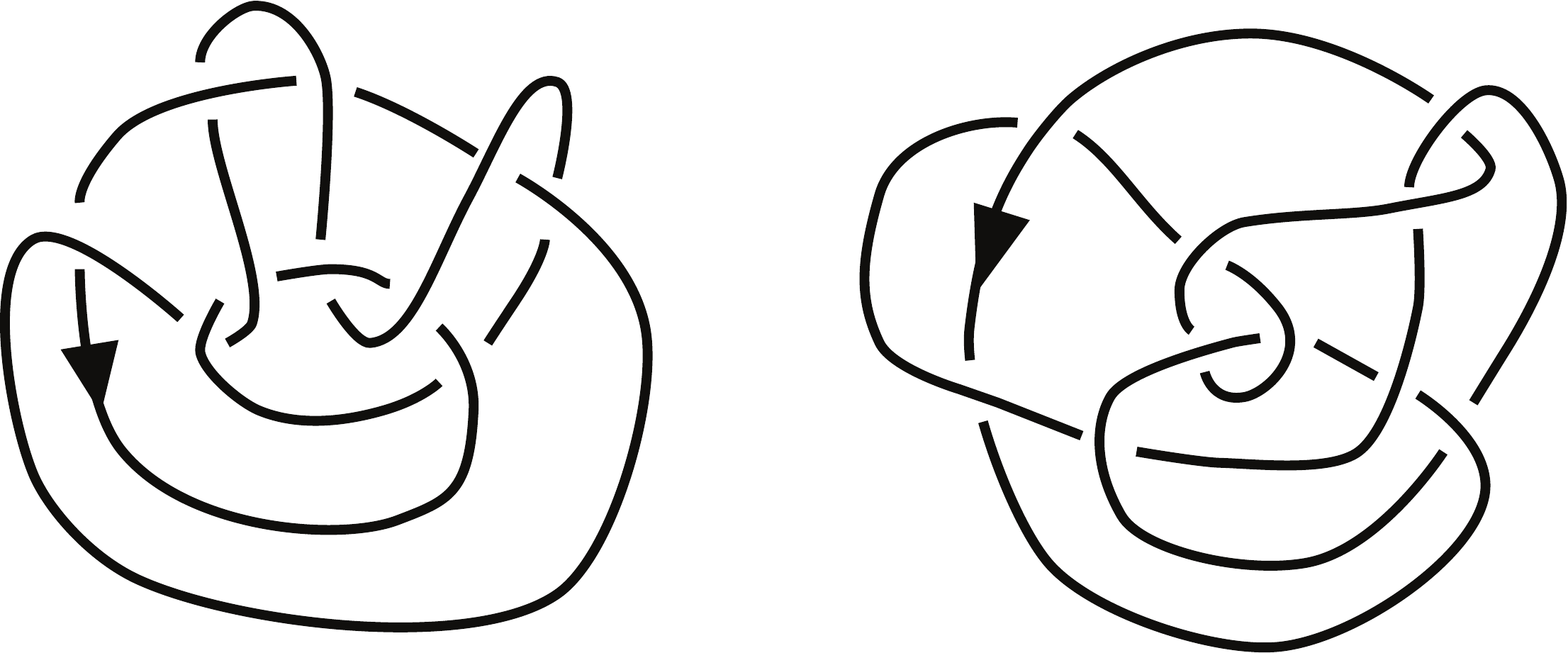}
\caption{Inserting tangle (a) into the precrossings in Figure~\ref{perko}.}
\label{perko_pos}
\end{center}
\end{figure}

Let us use the Alexander polynomial and take our tangle set $T$ to be the set consisting of tangles (a) and (b). If we insert tangle (b) at both precrossings in both diagrams, we recover the Perko pair itself, which were shown by Perko to be equivalent. On the other hand, the knot  $7_3$ is obtained by inserting tangle (a) at one precrossing and tangle (b) at the other precrossing (regardless of the order of insertion) in both diagrams. Finally, if tangle (a) is inserted at both precrossings in both diagrams as in Figure~\ref{perko_pos}, the two resulting knots are distinct. One knot is determined by the Alexander polynomial to be the knot $10_{148}$, while the other is knot $10_{160}$. Hence, the invariant $\overline{\mathcal{I}}$ derived from the Alexander polynomial and tangle set $T$ (where, say, $A_1=A_2=\frac{1}{2}$) can distinguish these two pseudoknots. \\

\noindent {\bf Remark.}  Notice that, because our pseudoknot invariant was defined using the Alexander polynomial and tangles (a) and (b) with the coefficient choice $A_1=A_2=\frac{1}{2}$, the were-set together with the Alexander polynomial can alternatively be used to distinguish the Perko pseudoknots.  \\

\begin{figure}[htbp]
\begin{center}
\includegraphics[height=1.2in]{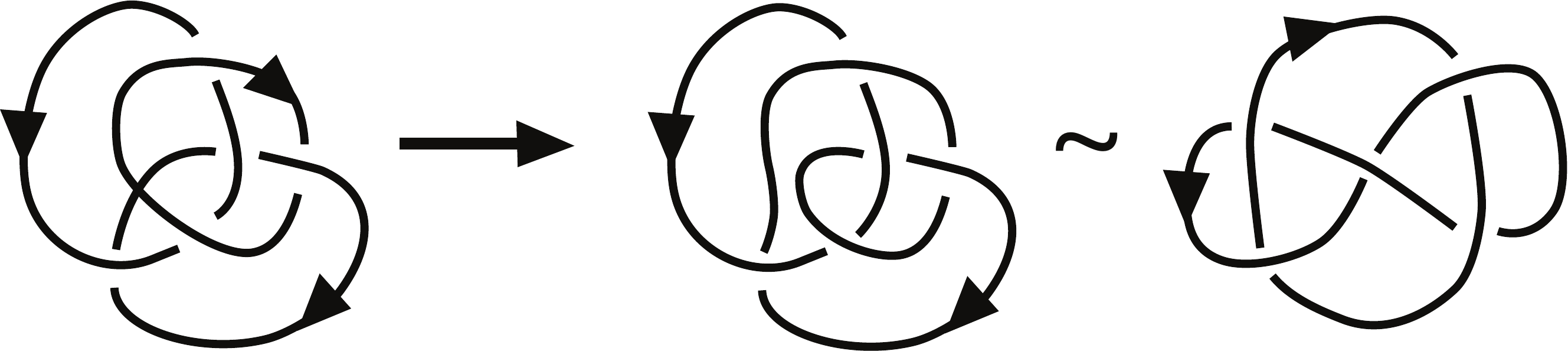}
\caption{Insertion of tangle (c) into a pseudo-Borromean ring diagram yields the Whitehead link, link $L5a1\{0\}$.}
\label{borr_insert}
\end{center}
\end{figure}

In our next example, we consider the pseudodiagram related to the Borromean rings pictured in Figure~\ref{borr_insert} (on the left). Here, we use the tangle insertion invariant $\overline{\mathcal{I}}$ that is defined using the Jones polynomial (computed via the bracket, with variable $A$) and take our tangle set to be the singleton set containing the tangle (c). As shown in Figure~\ref{borr_insert}, inserting the tangle (c) produces the Whitehead link. Since this link has Jones polynomial $$-A^6+A^{2}-2A^{-2}+A^{-6}-2A^{-10}+A^{-14},$$ $\overline{\mathcal{I}}$ can be used to prove that the original Borromean pseudoknot is nontrivial. Note that, if we had used a linking number invariant rather than the Jones polynomial to define our pseudoknot invariant, we would not have been able to detect the nontriviality of this example since the linking number of the Whitehead link is 0.

These are just a few examples to illustrate how invariants of pseudoknots can be derived from classical link invariants and tangle sets. We intend to consider many more examples in future work.

\section{Questions \& Future Work}

Given this framework for developing pseudoknot invariants, there are many questions that have yet to be be explored. We provide the reader with an initial list of open questions.

\begin{enumerate}
\item Are there examples of pseudoknot pairs that can only be distinguished using a tangle set $T$ that contains more complex tangles than the basic tangles (a), (b), and (c) in Figure~\ref{example_t}? For instance, can tangles (d), (e), or (f) be used to construct more powerful pseudoknot invariants? 
\item Given two arbitrary distinct pseudoknots $K_1$ and $K_2$, does there exist a classical link invariant $\mathcal{I}$ and tangle set $T$ such that the corresponding pseudoknot invariant $\overline{\mathcal{I}}$ distinguishes $K_1$ and $K_2$?
\item Can our invariant schema be generalized to include invariants of other sorts, e.g. the matrix of linking numbers? (This can be done for the case of a single precrossing.)
\item Is there a relationship between the Gauss-diagrammatic pseudoknot invariants defined in~\cite{gauss} and the invariant schema presented here?
\item Can we determine if a given reduced pseudodiagram has the fewest number of precrossings among all equivalent diagrams? Specifically, can our tangle insertion invariants be used for this purpose?
\end{enumerate}

\end{document}